\documentclass[12pt,oneside]{amsart}

\usepackage{amssymb}  
\usepackage{latexsym} 
\usepackage{comment}
\usepackage[all]{xy}

\DeclareFontEncoding{OT2}{}{} 
\newcommand{\textcyr}[1]{%
 {\fontencoding{OT2}\fontfamily{cmr}\fontseries{m}\fontshape{n}\selectfont #1}}

\newcommand{\Sha}{{\mbox{\textcyr{Sh}}}}

\newcommand{\Z}{{\mathbb Z}}
\newcommand{\Q}{{\mathbb Q}}
\newcommand{\R}{{\mathbb R}}

\newcommand{\F}{{\mathbb F}}
\newcommand{\BP}{{\mathbb P}}

\newcommand{\To}{\longrightarrow}

\newcommand{\tensor}{\otimes}
\newcommand{\disc}{\operatorname{disc}}
\newcommand{\Sel}{\operatorname{Sel}}
\newcommand{\Pic}{\operatorname{Pic}}

\newcommand{\fake}{{\text{fake}}}

\newcommand{\rank}{\operatorname{rank}}

\newcommand{\mult}{*}

                      {\hspace*{\fill}\nobreak$\Box$\par}

   {\begin{list}{}{\settowidth{\labelwidth}{#1}%
                   \setlength{\leftmargin}{\labelwidth}%
                   \addtolength{\leftmargin}{\labelsep}}}%
   {\end{list}}

\newtheorem{Theorem}{Theorem}[section]

\theoremstyle{definition}

\newtheorem{Question}[Theorem]{Question}

\numberwithin{equation}{section}

\newcounter{nootje}
\renewcommand\check[1]
  {\marginpar{\tiny\begin{minipage}{20mm}\begin{flushleft}\thenootje : #1\end{flushleft}\end{minipage}}\addtocounter{nootje}{1}}
\begin{document}

\title[Deciding existence of rational points]%
      {Deciding existence of rational points on curves: an experiment}

\author{Nils Bruin}
\address{Department of Mathematics,
         Simon Fraser University,
         Burnaby, BC,
         Canada V5A 1S6}
\email{nbruin@cecm.sfu.ca}
\thanks{Research of the first author supported by NSERC}

\author{Michael Stoll}
\address{School of Engineering and Science,
         International University Bremen,
         P.O.Box 750561,
	 28725 Bremen, Germany.}
\email{m.stoll@iu-bremen.de}
\date{April 25, 2006}


\maketitle


\section{Introduction}

The problem to decide whether a given algebraic variety defined over
the rational numbers has rational points is fundamental in Arithmetic
Geometry. Abstracting from concrete examples, this leads to the question
whether there exists an algorithm that is able to perform this task for any 
given variety. This is probably too much to ask for: we know that Hilbert's
Tenth Problem, which asks the same question for integral points on general
affine varieties, has a negative answer. But we can hope for a more
favorable outcome if we restrict the class of varieties we consider.

It is then most natural to look at curves first, since they have been
studied very intensely, resulting in a good theoretical knowledge and
a very rich supply of algorithmic methods. Also, it makes sense to consider 
the geometrically
nicest class of varieties, namely those that are projective. Since it is easy 
to check whether a curve has rational singular points, we can assume that
the curve is smooth.
Therefore, the question we are specifically interested in is the following.

\begin{Question}
  Is there an algorithm that decides for any given smooth projective curve
  $C/\Q$ whether $C$ has rational points or not?
\end{Question}

Since we can always algorithmically prove that $C(\Q) \neq \emptyset$ if
rational points exist by simply enumerating all rational points of the
relevant projective space and checking for each point if it is on~$C$, until
we find a rational point on~$C$, our question is equivalent to the
following, seemingly more restricted version.

\begin{Question}
  Is there an algorithm that verifies that $C(\Q) = \emptyset$ for any
  given smooth projective curve $C/\Q$ without rational points?
\end{Question}

``Verification'' here means that the algorithm constructs a proof of some
kind.

For curves of genus~$0$, our question has a positive answer since the
Hasse Principle holds for these curves: a curve of genus~$0$ has rational
points if and only if it is ``everywhere locally solvable'' (ELS), i.e., it
has real points and $p$-adic points for all primes~$p$.
Since for a general curve~$C$, we can check algorithmically whether it
has points everywhere locally, we can assume that $C$ is ELS. The
main problem is then to show that $C(\Q)$ is empty even though $C$ is ELS.

If $C$ is a curve of genus~$1$ with
Jacobian elliptic curve $E$, then we can perform descent calculations
(on~$E$ or on~$C$), which will succeed in proving that $C(\Q)$ is empty
if $C$ represents an element of~$\Sha(E)$ that is not divisible. In 
particular, if we assume that $\Sha(E)$ is finite for all elliptic
curves $E/\Q$, then our question has a positive answer for curves of genus~$1$
as well.

We will therefore focus our attention to curves of higher genus. 
It is only since fairly recently that there is some confidence that the
question might have a positive answer, spurred by progress on the theoretical
side~\cite{StollCov,StollOw} and also by heuristcal 
considerations~\cite{PoonenHeur}. In this paper, we attempt to give
supporting evidence of a more practical kind, by applying the available
algorithms (with some new improvements and additions) to a large number
of curves in order to see if we actually can decide for each of them
whether they have rational points or not.

The obvious class of curves to look at for a first attempt at gathering
evidence are the curves of genus~$2$. Their main advantage is that
quite a variety of algorithms are available for them, and so we can
hope to use them as adequate test cases. In order to keep the computational
effort within reasonable limits, we decided to consider ``small'' genus~$2$
curves. More precisely, our initinal set of curves consists of all
genus~$2$ curves over~$\Q$ that have a model of the form
\[ y^2 = f(x) =  f_6 x^6 + f_5 x^5 + f_4 x^4 + f_3 x^3 + f_2 x^2 + f_1 x + f_0 
\]
with integral coefficients $f_0, f_1, \dots, f_6$ satisfying $|f_j| \le 3$.
Excluding non-squarefree $f$ or $f$ of degree $\le 4$ and identifying
isomorphic curves, our initial set contains $196\,211$
isomorphism classes of curves.

In Section~\ref{sec:results}, we describe our findings, and in
Section~\ref{sec:methods},
we give an overview of the methods we have used. The details on the new methods
and the improvements on existing methods we have made can be found in a series
of forthcoming papers~\cite{BruinStoll2D,BruinStollFG,BruinStollBM}.

A complete list of all curves considered and lists of indications on
how to prove that each curve does or does not have rational points are available
at \cite{BruinStoll:data}.

\subsection*{Acknowledgments}

We would like to thank Victor Flynn and Bjorn Poonen
for useful discussions related to our project. M.S. thanks the
Computational Laboratory for Analysis, Modeling and Visualization (CLAMV) of
International University Bremen for the possibility to use computing time
on the CLAMV Teaching Lab machines. This was used for substantial parts of 
the computations that were done in the course of this project. 
For the computations, the {\sf MAGMA}~\cite{Magma} system was used.


\section{Results}\label{sec:results}

In a first step, we searched for a small rational point on each curve~$C$.
Note that $C$ has one or two obvious points if $f_0 \in \{0, 1\}$
or $f_6 \in \{0, 1\}$ ($C$ is considered to have one or two rational points
``at infinity'' if $f_6 = 0$ or a non-zero square, respectively).
At a later stage, we searched for larger rational points on those curves
that were not yet decided. The largest points found at this stage were
$(1519/601, 4816728814/601^3)$ on
\[ C : y^2 = 3\,x^6 - 2\,x^5 - 2\,x^4 - x^2 + 3\,x - 3 \]
and $(193/436, 165847285/436^3)$ on
\[ C : y^2 = 3\,x^6 - 3\,x^5 - x^4 - x^3 - 3\,x^2 + x - 3 \,. \]

This left us with $58\,681$ curves~$C$ without (apparent) rational points, for
which we need to prove that $C(\Q) = \emptyset$. Among these, there are
$29\,278$ curves with points everywhere locally. 
Together with the curves that do have rational points,
this means we found $166\,808$ with points everywhere locally, which is about
$85$\,\% of all the curves we considered. In~\cite{PoonenStoll2} 
(see also~\cite{PoonenStoll1}, Section~9), it is shown that the set of
polynomials~$f$ giving rise to an everywhere locally solvable curve has a 
well-defined, positive density~$\delta$. Numerical estimates of the local 
densities involved lead to a value close to $0.85$ for~$\delta$, which fits
well with our observations.

The next stage in the procedure is to perform a $2$-cover descent on
each of the remaining curves. This constructs (implicitly) a finite collection
of curves $D_j$ that cover~$C$ and such that every rational point on~$C$
is the image of a rational point on some~$D_j$. So if we obtain an empty
covering collection $\{D_j\}$, this proves that $C$ has no rational points.
For a more precise description of the computation, see Section~\ref{sec:2cov}.
With this method, we were able to prove that $C(\Q)$ is empty for all but
$1492$~curves.

For these $1492$~curves, we wanted to perform a ``Mordell-Weil sieve''
computation. The idea is as follows. Let $J$ be the Jacobian variety of~$C$,
and assume that we can embed $C$ into~$J$. Assume also that we can determine
generators of the Mordell-Weil group~$J(\Q)$. Now consider the following
commutative diagram.
\[ \xymatrix{ C(\Q) \ar[r] \ar[d] & J(\Q) \ar[d]^-{\alpha} \\
              \prod\limits_{p \in S} C(\F_p) \ar@<4pt>[r]^{\beta} 
                & \prod\limits_{p \in S} J(\F_p)
            }
\]
Here $S$ is some finite set of primes (of good reduction for~$C$, say).
Since we know~$J(\Q)$ and can find the finite sets $C(\F_p)$, we can 
compute the images of $\alpha$ and~$\beta$. If these images do not meet,
this proves that $C(\Q)$ is empty.

As a first step, we had to find generators of the Mordell-Weil group.
To do this, we first perform a 2-descent (see~\cite{Stoll2DH}) on~$J$
to get an upper bound for the rank of the finitely generated abelian
group~$J(\Q)$. Then we need to find the correct number of independent
points in~$J(\Q)$. In order to be able to do this successfully, we had
to come up with new strategies, involving a search for points on (quotients) of
$2$-covering spaces for~$J$. See Section~\ref{sec:gens} for more details.
In this way, we were able to find generators of a finite index subgroup
of~$J(\Q)$ for all but 47 curves. It is then a fairly easy matter to
check that we actually had generators of~$J(\Q)$ (modulo torsion),
see~\cite{StollHt2}. In the course of these computations, we also
found a rational point on the 2-covering space $\Pic^1_C$ for~$J$,
which provides an embedding of $C$ into~$J$. So for these 1445~curves
(3 of rank~$0$, 516 of rank~$1$, 772 of rank~$2$, 152 of rank~$3$, 
and 2 of rank~$4$),
the assumptions for the application of the Mordell-Weil sieve are
satisfied. After several improvements of the algorithm performing the
actual sieve computation (the problem here is combinatorial explosion),
we were finally able to run the procedure successfully for all these
curves. With the current version of the algorithm, the maximal computation
time for a single curve was roughly 16~hours on a 1.7~GHz processor;
this curve is one of the two with Mordell-Weil rank~$4$.

For the remaining 47~curves (36 of Selmer rank~$2$, 10 of Selmer rank~$3$,
and one of Selmer rank~$4$), the number of independent points we found
fell short by~$2$ of the Selmer rank. Therefore we suspect that there
is nontrivial $2$-torsion in~$\Sha(J)$ in these cases. In 5 out of the 10
Selmer rank~$3$ cases, we found a rational point on~$\Pic^1_C$. Here
we expect that $\Sha(J)[2] = (\Z/2\Z)^2$. This was confirmed by an ad-hoc
visualization argument. See~\cite{B:vissha,BF:vissha} for a description and
detailed analysis of this method for hyperelliptic curves with a rational
branch point. 

For these 5 curves, we know that $\rank J(\Q) = 1$, and we have an embedding of
$C$ into~$J$, so that we can run the Mordell-Weil sieve procedure, which
confirms that there are no rational points on these curves.

One example we found is
\[ C : y^2 = f(x) = -x^6 + 2x^5 + 3x^4 + 2x^3 - x - 3 \,. \]
When we consider a quadratic twist of this curve,
\[ C^{(-1)} : y^2 = -f(x)\,, \]
we find that $J^{(-1)}(\Q)$ is of rank~$4$, where $J^{(-1)}$ is the
Jacobian of~$C^{(-1)}$.
A slightly more involved computation
gives that $J(\Q(\sqrt{-1}))$ is of rank at most~$5$. Since this rank is
the sum of the ranks of $J(\Q)$ and $J^{(-1)}(\Q)$, this means
that the rank of $J(\Q)$ can be at most~$1$. This is less than the rank
bound~$3$ we obtain from a $2$-descent on $J$ directly.

In the remaining 42~cases, we did not find a rational point on~$\Pic^1_C$.
On the other hand, from the $2$-cover descent, we know that $C$ has
everywhere locally solvable $2$-coverings; the same must then be true
for~$\Pic^1_C$. This means that the class of~$\Pic^1_C$ in~$\Sha(J)$ is
divisible by~$2$. If $\Pic^1_C(\Q) = \emptyset$, then this implies that
there are elements of order~$4$ in~$\Sha(J)$. The computations necessary for a
visualization argument are hardly feasible in this situation: one needs to
compute the $2$-Selmer group of~$J$ over a quartic number field. This
involves finding an $S$-unit group in a degree~$24$ number field.

Still, assuming GRH, we were successful for $4$~curves in showing that the
true Mordell-Weil rank is smaller than the bound obtained from a $2$-descent.
One of these curves is
\[ C : y^2 = -3x^6 - x^5 + 2x^4 + 2x^2 - 3x - 3\,. \]
The Jacobians of the quadratic twists by $2,-3,-6$ can easily be shown to have
Mordell-Weil ranks $4,4,3$ respectively. Furthermore, a $2$-descent
shows, conditional on GRH, that $J(\Q(\sqrt{2},\sqrt{-3}))$ is of rank at
most~$11$. It follows that $J(\Q)$ must be of rank~$0$.

We do not expect that results along these lines can be extended much further. To
complement the above computations, we computed the analytic order of~$\Sha(J)$.
For this we had to assume that the $L$-series $L(C,s)$ can be analytically
continued and satisfies the usual functional equation. The results
of our computations are consistent with this assumption. First of all, we
verified that the $r$-th derivative of $L(C,s)$ at~$s=1$ is nonzero, where
$r$ is the conjectured rank of~$J(\Q)$ (i.e., the number of independent
points we have found). Secondly, the analytic
order of~$\Sha(J)$ comes out to be~$16$ for the 42~curves where we
expect elements of order~$4$, and it is~$4$ for the 5~curves mentioned
above, where we expect $\Sha(J) = (\Z/2\Z)^2$. Hence, assuming standard
conjectures on $L$-series and the Birch and Swinnerton-Dyer Conjecture,
we find that $\Pic^1_C(\Q) = \emptyset$ for our 42~curves, and therefore
$C(\Q) = \emptyset$ as well.

\medskip

The main result of our experiment is that
we were successful in deciding the existence of rational
points unconditionally on all but~42 of our curves.
If we assume standard conjectures, we can prove that there are no rational
points on these remaining 42~curves as well.
This very positive result lends strong support to the belief that
existence of rational points on curves should be decidable,
at least for curves of genus $2$.

Note that our results also provide evidence for the conjecture that
the Brauer-Manin obstruction should be the only obstruction against
rational points on curves. For all but the 1492 curves surviving a
$2$-cover descent, we verify this unconditionally. For the remaining
curves, we need to assume that $\Sha(J)$ has trivial
divisible subgroup, plus whatever assumptions were necessary in addition
for individual curves. See~\cite{Scharaschkin,Flynn,StollCov} for 
details on how our computations relate to the Brauer-Manin obstruction.

\renewcommand{\arraystretch}{1.2}

\begin{table}
  \begin{tabular}{|l|r|r|}
    \hline
    All curves                             & 196\,211 & 100.00\,\% \\
    Curves with rational points            & 137\,530 &  70.09\,\% \\
    Curves without rational points         &  58\,681 &  29.91\,\% \\
    ELS curves total                       & 166\,808 &  85.01\,\% \\
    ELS curves without rational points     &  29\,278 &  14.92\,\% \\
    Curves with ELS $2$-covers among these &   1\,492 &   0.76\,\% \\
    Curves that need GRH or BSD conjecture &       42 &   0.00\,\% \\
    \hline
  \end{tabular}

\bigskip
\caption{Curve statistics (ELS = everywhere locally solvable)}

\[ \begin{array}{|l||r|r|r||r|}
     \hline 
     \text{conj.\ } \Sha(J) &   0 & (\Z/2\Z)^2 & (\Z/4\Z)^2 & \text{Total} \\
     \hline
     \rank J(\Q) = 0        &   3 &            &         36 &   39 \\
     \rank J(\Q) = 1        & 516 &          5 &          5 &  526 \\
     \rank J(\Q) = 2        & 772 &            &          1 &  773 \\
     \rank J(\Q) = 3        & 152 &            &            &  152 \\
     \rank J(\Q) = 4        &   2 &            &            &    2 \\
     \hline
     \text{all ranks}      & 1445 &          5 &         42 & 1492 \\
     \hline
  \end{array}
\]
\smallskip
\caption{Ranks and conjectural $\Sha$ for the curves surviving 
         $2$-cover descent}
\end{table}

\section{Methods}\label{sec:methods}

In this section, we give an overview of the methods we have used. 
Detailed descriptions will be provided 
in~\cite{BruinStoll2D,BruinStollFG,BruinStollBM}.


\subsection{$2$-Cover Descent}\label{sec:2cov}

Let the curve~$C$ be given by the equation $y^2 = f(x)$, and let
$L$ denote the \'etale $\Q$-algebra $\Q[T]/(f(T))$. We let $\theta$
be the image of~$T$ in~$L$. If $f$ has a rational root or is of odd degree
then $C$ has a rational point. Therefore, we can assume in the following that
$f$ is of degree $6$ and has no rational roots.
Let $a$ be the leading coefficient
of~$f$. Let $k = \Q$ or $k = \Q_v$, where $v$ is a prime~$p$ or~$\infty$
and $\Q_\infty = \R$. We have a map
\[ F : C(k) \To \frac{(L \tensor_{\Q} k)^\mult}%
                     {k^\mult (L \tensor_\Q k)^{\mult2}} \,, \qquad
       P = (x,y) \longmapsto
          (x - \theta) \cdot k^\mult (L \tensor_\Q k)^{\mult2} \,,
\]
whose image is contained in the subset of elements whose norm
in~$k^\mult/k^{\mult2}$ is the class of~$a$. Note that
$y^2 = f(x) = a N_{L \tensor k/k}(x - \theta)$.

As in the case of $2$-descent
on the Jacobian~$J$ of~$C$, one shows that $F(C(\Q_p))$ is contained
in the image of the $p$-adic units for all odd~$p$ not dividing the
discriminant of~$f$ (see~\cite{Stoll2DH}). Let 
$H' \subset L^\mult/\Q^\mult L^{\mult2}$ be the (finite) set of
elements that come from $p$-adic units for this set of primes, and 
let $H \subset H'$ be the subset of elements whose
norm is $a \Q^{\mult2}$; then $F(C(\Q)) \subset H$. There are cases
when $H$ is already the empty set; we can then immediately conclude 
that $C(\Q) = \emptyset$. An example of this~is
\[ C : y^2 = 2x^6 + 3x^5 + x^4 - 3x^3 - 2x^2 + 2x + 3\,. \]
For this curve it can be checked that $2 \Q^{\mult2}$ is not the norm
of an element of~$H'$.

We denote $L \tensor_{\Q} \Q_v$ by~$L_v$.
Let $H_v \subset L_v^\mult/\Q_v^\mult L_v^{\mult2}$ denote the subset
of elements whose norm is $a \Q_v^{\mult2}$. We have the following
commutative diagram.
\[ \xymatrix{ C(\Q) \ar[r]^{F} \ar[d] & H \ar[d]^-{\rho} \\
              \prod\limits_{v \in S} C(\Q_v) \ar@<4pt>[r]^{F} 
                 & \prod\limits_{v \in S} H_v
            }
\]
Here, $S$ is a suitable finite set of places. One can show that 
$F(C(\Q_p)) = H_p$ for $p > 1154$ if $p$ does not divide~$\disc(f)$.
Therefore, we obtain the maximal amount of information when we choose
\[ S = \{\infty\} \cup \{p : p < 1154\} \cup \{p : p \mid \disc(f)\} \,. \]
Note that the sets $H_v$ are finite and that $F$ is $v$-adically continuous,
hence locally constant. Therefore we can compute $F(C(\Q_v)) \subset H_v$
explicitly for every~$v$. Following \cite{PoonenSchaefer}, we define the 
{\em fake $2$-Selmer set} of~$C$, $\Sel_\fake^{(2)}(C)$,
to be the preimage in~$H$ under~$\rho$ of the image of the lower $F$ map.
Then $F$ maps $C(\Q)$ into~$\Sel_\fake^{(2)}(C)$, hence
if $\Sel_\fake^{(2)}(C) = \emptyset$, then we know that $C(\Q)$ is empty 
as well.

The geometric interpretation of the elements of~$\Sel_\fake^{(2)}(C)$ is that
they correspond to everywhere locally solvable $2$-covering curves of~$C$.
If $\xi \in L$ represents an element of~$\Sel_\fake^{(2)}(C)$, then the 
corresponding covering $D_\xi \to C$ can be obtained as follows.
We write $z = z_0 + z_1 \theta + \dots + z_5 \theta^5$ for a generic
element of~$L$. The condition for a rational point~$P = (x,y)$ on~$C$ to be
in the image of~$D_\xi(\Q)$ is that 
\[ (x - \theta) \cdot \Q^\mult L^{\mult2}
     = F(P) = \xi \cdot \Q^\mult L^{\mult2} \,. 
\]
So $x - \theta = c \xi z^2$ for some $c \in \Q$, $z \in L$. Expanding the right
hand side in terms of powers of~$\theta$, we obtain four quadrics in the
six variables $z_0, \dots, z_5$ that express the condition that the 
coefficients of $\theta^2, \dots, \theta^5$ have to vanish. These four
quadrics define the curve~$D_\xi \subset \BP^5$ of degree~$16$ 
and genus~$17$. To obtain the
covering map, note that $x$ can be recovered from the coefficients of
$1$ and~$\theta$ in~$\xi z^2$, and $y$ can be recovered from these, the
norm of~$z$ and a square root of $N(\xi)/a$.
One has to make a sign choice here, so that there are really
two different covering maps in most cases. See also \cite[5.3]{Bruin:Thesis} and
\cite{BruinFlynn:Tow2Cov} for a description of the cover.
For details on how to compute $\Sel_\fake^{(2)}(C)$ efficiently,
see~\cite{BruinStoll2D}.


\subsection{Finding Generators}\label{sec:gens}

Since the simplest generally available model of the Jacobian~$J$ is
given by 72~quadrics in~$\BP^{15}$ (see~\cite{CasselsFlynn}), 
it is usually not a good idea to
search for rational points directly on~$J$. A better alternative is
to consider the Kummer surface $K = J/\{\pm 1\}$, which sits naturally
as a quartic surface in~$\BP^3$. We now can search for rational points
on~$K$ that lift to rational points on~$J$. A fairly efficient implementation
of this idea that uses mod-$p$ information for several primes~$p$ in
order to rule out many candidates is obtainable as {\tt j-points}
from M.~Stoll's homepage; this program is also incorporated in~{\sf MAGMA}.
This approach is feasible for points of naive (non-logarithmic) height
around~10\,000 or a little bit more. (The height is that of the image
point on~$K \subset \BP^3$.)

However, there are many cases in our list where there is a much bigger
generator. In order to find these, we use the idea (by now in common use
in the context of elliptic curves) that rational points on~$J$ lift
to rational points on a $2$-covering of~$J$ that should be much smaller.
Therefore we attempt to search for rational points on these $2$-coverings.
However, these coverings are as complicated geometrically as $J$ itself,
therefore we consider a suitable quotient again. 

Recall (see~\cite{PoonenSchaefer,Stoll2DH}) that the fake $2$-Selmer group 
of~$J$
is a finite subgroup of $L^\mult/\Q^\mult L^{\mult2}$. It contains
the image of~$J(\Q)$ under a map that sends a rational point~$P$ to an
element represented by $x_0 - x_1 \theta + \theta^2 \in L$, for
certain $x_0, x_1 \in \Q$ depending on~$P$. Let $\xi$ be an element
of the fake Selmer group. We use the same idea as in the previous section
to construct a surface~$K_\xi$: we are looking for $z \in L$ such that
$\xi z^2$ does not involve $\theta^3, \theta^4, \theta^5$. This gives
us an intersection~$K_\xi$ of three quadrics in~$\BP^5$. We simplify the 
defining
equations as far as possible by a change of projective coordinates so that 
they have small coefficients. Then we perform a search for rational
points on~$K_\xi$ using a $p$-adic variant of Elkies' lattice-based point
searching techniques (see \cite{Elkies}). For each point found, we check 
whether it corresponds to a rational point on~$J$. In this way, we can find 
points in~$J(\Q)$ whose image in the fake Selmer group is nontrivial.

However, note that $\Pic^1_C$ is a $2$-covering
of~$J$ via the map $D \mapsto 2D-W$, where $D\in\Pic^1_C$ and $W$ is the
canonical class. Its image in the fake Selmer group is trivial, so the method
above will not help in finding rational points on it.
Instead, in analogy to the use of the Kummer surface when searching for points
on~$J$, we can use the
dual Kummer surface (see~\cite{CasselsFlynn}, Chapter~4). We can even
go a step further and consider $2$-coverings of~$\Pic^1_C$. In this
case, we obtain 3-dimensional varieties, given as intersections
of two quadrics in~$\BP^5$, that are quotients of $\BP^1$-bundles over
the coverings we are interested in. We can search for rational points
on these 3-folds and check whether they give rise to a rational point
on~$J$. This amounts to a partial explicit $4$-descent on~$J$. It is
therefore perhaps not surprising that we were able to find some quite
large generators in this way. The record example is
\[ C : y^2 = -3\,x^6 + x^5 - 2\,x^4 - 2\,x^2 + 2\,x + 3 \]
with $J(\Q)$ infinite cyclic generated by $P_1 + P_2 - W$ where the
$x$-coordinates of $P_1$ and~$P_2$ are the roots of
\[ x^2 + \tfrac{37482925498065820078878366248457300623}%
               {34011049811816647384141492487717524243}\,x 
       + \tfrac{581452628280824306698926561618393967033}%
               {544176796989066358146263879803480387888} \,;
\]
the canonical logarithmic height of this generator is~$95.26287$.
The second largest example is
\[ C : y^2 = -2\,x^6 - 3\,x^5 + x^4 + 3\,x^3 + 3\,x^2 + 3\,x - 3 \]
with $J(\Q)$ generated by a point coming from
\[ x^2 + \tfrac{83628354341362562860799153063}{26779811954352295849143614059}\,x
       + \tfrac{852972547276507286513269157689}{321357743452227550189723368708}
       \,.
\]
The canonical height of this generator is~$77.33265$.

For details see~\cite{BruinStollFG}.


\subsection{Mordell-Weil Sieve}\label{sec:sieve}

As mentioned in Section~\ref{sec:results}, we consider the commutative diagram
\[ \xymatrix{ C(\Q) \ar[r] \ar[d] & J(\Q) \ar[d]^-{\alpha} \\
              \prod\limits_{p \in S} C(\F_p) \ar@<4pt>[r]^{\beta} 
                & \prod\limits_{p \in S} J(\F_p)
            }
\]
with a suitable finite set~$S$ of (good) primes. In some cases, it can
be helpful to use some more general finite quotient of~$J(\Q_p)$ instead
of~$J(\F_p)$, for example to make use of information modulo higher powers
of~$p$, or also in order to use information at primes of bad reduction.
In the following discussion, we will assume for simplicity that we
are working with~$J(\F_p)$.

Our goal is to prove that the images of $\alpha$ and~$\beta$ above do not
meet, for some set~$S$, which implies that $C(\Q) = \emptyset$. This approach
was (to our knowledge) first suggested by Scharaschkin~\cite{Scharaschkin}.
Flynn~\cite{Flynn} used it for more extensive calculations. We
would like to mention here that in the course of improving the algorithms,
we were able to prove all the curves marked ``Unresolved'' in the tables
of~\cite{Flynn} to have no rational points. (All but five of these already
succumb to a 2-cover descent, the remaining five, all of which have
Jacobians of Mordell-Weil rank~$3$, can be dealt with using our
Mordell-Weil sieve implementation.)

The basic algorithmic problem one has in this computation is that
the product of the~$J(\F_p)$ can be a very large group. One approach
to keep the combinatorics in check is to work with $J(\Q)/B J(\Q)$
and $\prod_p J(\F_p)/B J(\F_p)$ for a suitable choice of~$B$. In practice,
we compute the subset of~$J(\Q)/B J(\Q)$ that maps under~$\alpha$ into
the image of~$\beta$. We first need to choose a promising set~$S$ of
primes. Since we can only hope to arrive at a contradiction when the
group orders of the~$J(\F_p)$ have (preferably large) common factors,
we select those primes~$p$ for which the order of~$J(\F_p)$ is sufficiently
smooth. We then compute the image of~$C(\F_p)$ in~$J(\F_p)$ and the
image of the generators of~$J(\Q)$. Note that this involves a discrete
logarithm computation in~$J(\F_p)$ for each point in~$C(\F_p)$ and each 
generator of~$J(\Q)$. While this is a hard problem in general,
it is harmless here since the group order is smooth and we can reduce
to several discrete logs in small groups.

In the following discussion, the set~$S$ is fixed.
For a given~$B$, we can find the image $C_{B,p}$ of~$C(\F_p)$ 
in~$J(\F_p)/B J(\F_p)$, and we can then compute the expected size
\[ n(B) = \#\bigl(J(\Q)/B J(\Q)\bigr)
            \prod_{p \in S} \frac{\# C_{B,p}}{\#\bigl(J(\F_p)/B J(\F_p)\bigr)}
\]
of the subset~$A(B)$ of~$J(\Q)/B J(\Q)$ that maps into these images for all 
$p \in S$.
We now search for a sequence $1 = B_0, B_1, \dots, B_m$ such that
$B_{j+1} = B_j q_j$ for some prime~$q_j$, such that~$n(B_m) \ll 1$
and such that $\max_j n(B_j)$ is not too large. (See~\cite{PoonenHeur}
for heuristics why there should exist $B$ with $n(B) \ll 1$, at least
when $S$ is sufficiently large.)

After we have fixed our sequence $(B_j)$, we successively compute the
sets~$A(B_j)$ for $j = 1, 2, \dots$ until $A(B_j) = \emptyset$. If we
reach $j = m$ and $A(B_m) \neq \emptyset$, then we can check if this is
caused by a exhibitable rational point.
The set $A(B_m)$ will give a very good indication of which elements 
of~$J(\Q)$ could
give rise to such a point. If we cannot find a point, we can extend the
sequence or choose a bigger set~$S$. In fact,
this situation never occurred in our computations.

To obtain $A(B_{j+1})$
from~$A(B_j)$, we run through the elements of~$A(B_j)$. For each element,
we run through its possible lifts to $J(\Q)/B_{j+1} J(\Q)$, and for each
lift check whether it maps into the image of~$C$ mod~$p$ for all relevant~$p$
(i.e., such that the largest power of~$q_{j+1}$ dividing $B_{j+1}$ also
divides the exponent of~$J(\F_p)$). The largest set $A(B_j)$ that we
encountered in our computations had a size of about~$10^6$. It is perhaps worth
while mentioning that the estimate~$n(B)$ for $\#A(B)$ was in most
cases accurate up to a factor of 2 to~5, so that a value 
$n(B_m) < 10^{-3}$ (say) virtually guarantees success in practice.

For details see~\cite{BruinStollBM}.


\subsection{BSD Computations}\label{sec:bsd}

Finally, let us give some indications of how to compute the analytic
order of~$\Sha$. Dokchitser~\cite{Dokchitser} describes how the numbers
$L^{(r)}(C, 1)$ can be computed numerically, given (i) the coefficients
$a_n$ of the $L$-series for sufficiently many~$n$, (ii) the conductor~$N$
of~$C$ (or~$J$), and (iii) the sign $\varepsilon$ in the (conjectured)
functional equation. The latter is determined by the parity of the rank.
The coefficients $a_p$ and $a_{p^2}$ for good primes~$p$ can be computed
by counting the points in~$C(\F_p)$ and~$C(\F_{p^2})$; these coefficients
then determine $a_{p^k}$ for all $k \ge 1$. For bad primes~$p$,
the coefficients can be deduced from a minimal proper regular model
of~$C$ over~$\Z_p$; a description of the computation of such a model
can be found in~\cite{6Authors}. The most frequent case is that an odd prime~$p$
divides the discriminant of the polynomial~$f$ just once; then
\[ f(x) \equiv (x - a)^2 g(x) \bmod p \,, \]
and the Euler factor at~$p$ of~$L(C, s)$ depends on whether $g(a)$ is
a square or not and the number of $\F_p$-points on the genus~$1$ curve
$y^2 = g(x)$. In most other cases, the original model is already regular.
For all of the curves, we computed $5 \cdot 10^5$ or even $10^6$ coefficients;
this led to an error in the value of $\#\Sha(J)$ predicted by the
Birch and Swinnerton-Dyer Conjecture of less than $10^{-3}$ in all cases.

We can find the odd part of the conductor~$N$ using Q.~Liu's 
{\tt genus2reduction} program~\cite{LiuProgram}, based on~\cite{Liu}.
If the given model of the curve is regular at~$2$, then the power of~$2$
dividing~$N$ is that dividing the discriminant of~$C$. Otherwise, we use
the approach described in~\cite[\S\,7]{Dokchitser} to determine the right
power of~$2$ (which is then less than that in the discriminant). We
can then verify the functional equation for the inverse Mellin transform
of~$L(C, s)$ numerically, thus corroborating our computations.

Given the value of $L^{(r)}(C, 1)$, we compute the analytic order of~$\Sha$
by solving the conjectural equality between $L^{(r)}(C, 1)/r!$ and
a combination of invariants of~$C$ and~$J$ for $\#\Sha$. 
See~\cite{6Authors} for how to compute the other invariants.
As already mentioned, the values we obtain were always close to an
integer, which was~$4$ in the five cases where we expect 
$\Sha(J) \cong (\Z/2\Z)^2$ and~$16$ in the remaining cases, where we
expect $\Sha(J) \cong (\Z/4\Z)^2$.



\begin{thebibliography}{MM99}
\frenchspacing
\bibitem[Br1]{Bruin:Thesis}
  {\sc N. Bruin:} 
  {\it Chabauty methods and covering techniques applied to generalized Fermat
  equations},
  CWI Tract 133, 77 pages (2002). 

\bibitem[Br2]{B:vissha}
  {\sc N. Bruin:}
  {\it Visualisation of Sha[2] in Abelian Surfaces,}
  Math. Comp. {\bf 73}, no. 247, 1459--1476 (2004).

\bibitem[BF1]{BruinFlynn:Tow2Cov}
  {\sc N. Bruin} and {\sc E.V. Flynn:}
  {\it Towers of 2-covers of hyperelliptic curves,}
  Trans. Amer. Math. Soc. {\bf 357}, 4329--4347 (2005).

\bibitem[BF2]{BF:vissha}
  {\sc N. Bruin} and {\sc E.V. Flynn:}
  {\it Exhibiting Sha[2] on Hyperelliptic Jacobians,}
  to appear in Journal of Number Theory (2005?).

\bibitem[BS1]{BruinStoll2D}
  {\sc N. Bruin} and {\sc M. Stoll:}
  {\it 2-cover descent on hyperelliptic curves,}
  in preparation.

\bibitem[BS2]{BruinStollFG}
  {\sc N. Bruin} and {\sc M. Stoll:}
  {\it Finding Mordell-Weil generators on genus 2 Jacobians,}
  in preparation.
  
\bibitem[BS3]{BruinStollBM}
  {\sc N. Bruin} and {\sc M. Stoll:}
  {\it The Mordell-Weil sieve: Proving non-existence of rational points on 
  curves,}
  in preparation.

\bibitem[BS4]{BruinStoll:data}
  {\sc N. Bruin} and {\sc M. Stoll:}
  Electronic data,
  available from \\
  \texttt{http://www.cecm.sfu.ca/\textasciitilde nbruin/smallgenus2curves}.
  
\bibitem[CF]{CasselsFlynn}
  {\sc J.W.S. Cassels} and {\sc E.V. Flynn:} 
  {\it Prolegomena to a middlebrow arithmetic of curves of genus~2,} 
  Cambridge University Press, Cambridge (1996).

\bibitem[Do]{Dokchitser}
  {\sc T. Dokchitser:}
  {\it Computing special values of motivic L-functions,}
  Exp. Math. {\bf 13}, No.2, 137--149 (2004).

\bibitem[El]{Elkies}
  {\sc N. D. Elkies:} {\it Rational points near curves and small nonzero
  $\vert x\sp 3-y\sp 2\vert $ via lattice reduction,} in:
  {\sc W. Bosma} (ed.): {\it Algorithmic number theory (Leiden, 2000),}
  Springer, Berlin, Lecture Notes in Comput. Sci. {\bf 1838}, 33--63 
  (2000).
  
\bibitem[Fl]{Flynn}
  {\sc E.V. Flynn:} {\it The Hasse Principle and the Brauer-Manin obstruction 
  for curves}, Manuscripta Math. {\bf 115}, 437--466 (2004).

\bibitem[F+]{6Authors}
  {\sc E.V. Flynn, F. Lepr\'evost, E.F. Schaefer, W.A. Stein, M. Stoll}
  and {\sc J.L. Wetherell:} {\it Empirical evidence for the Birch and 
  Swinnerton-Dyer conjectures for modular Jacobians of genus 2 curves},
  Math.\ Comp. {\bf 70}, 1675--1697 (2001).

\bibitem[Li1]{Liu}
  {\sc Q. Liu:} {\it Conducteur et discriminant minimal de courbes de genre~2,}
  Compositio Math. {\bf 94}, 51--79 (1994).

\bibitem[Li2]{LiuProgram}
  {\sc Q. Liu:} {\tt genus2reduction} program, available at\\
  \verb+http://www.math.u-bordeaux.fr/~liu/G2R/+.

\bibitem[M]{Magma} {\sf MAGMA} is described in
  {\sc W. Bosma, J. Cannon} and {\sc C. Playoust:} 
  {\it The Magma algebra system I: The user language}, 
  J. Symb. Comp. {\bf 24}, 235--265 (1997).
  (Also see the Magma home page at 
  {\tt http://www.maths.usyd.edu.au:8000/u/magma/}\,.)

\bibitem[Po]{PoonenHeur}
  {\sc B. Poonen:} 
  {\it Heuristics for the Brauer-Manin obstruction for curves,}
  Preprint (2005).

\bibitem[PSc]{PoonenSchaefer}
  {\sc B. Poonen} and {\sc E.F. Schaefer:} 
  {\it Explicit descent for Jacobians of cyclic covers of the projective line,} 
  J. reine angew. Math. {\bf 488}, 141--188 (1997).

\bibitem[PSt1]{PoonenStoll1}
   {\sc B. Poonen} and {\sc M. Stoll:} 
   {\it The Cassels-Tate pairing on principally polarized abelian varieties,}
   Ann. of Math. {\bf 150}, 1109--1149 (1999).

\bibitem[PSt2]{PoonenStoll2}
   {\sc B. Poonen} and {\sc M. Stoll:} 
   {\it A local-global principle for densities,}
   in: {\sc Scott D. Ahlgren} (ed.) et al.: 
   {\it Topics in number theory. In honor of B. Gordon and S. Chowla.}
   Kluwer Academic Publishers, Dordrecht. 
   Math. Appl., Dordr. {\bf 467}, 241--244 (1999).

\bibitem[Sc]{Scharaschkin}
  {\sc V. Scharaschkin:} {\it The Brauer-Manin obstruction for curves},
  Manuscript.

\bibitem[St1]{Stoll2DH}
  {\sc M. Stoll:} {\it Implementing 2-descent on Jacobians of hyperelliptic
  curves,} Acta Arith. {\bf 98}, 245--277 (2001).

\bibitem[St2]{StollHt2}
  {\sc M. Stoll:} {\it On the height constant for curves of genus two, II,}
  Acta Arith. {\bf 104}, 165--182 (2002).
  
\bibitem[St3]{StollCov}
  {\sc M. Stoll:} {\it Finite descent and rational points on curves,}
  Preprint (2006).

\bibitem[St4]{StollOw}
  {\sc M. Stoll:} {\it Finite coverings and rational points,} 
  Oberwolfach Report 32/2005 (2005).

\end{thebibliography}
\end{document}